\documentclass[11pt]{article}

\usepackage[margin=1.1in]{geometry}
\usepackage{amsmath,amssymb,amsthm,mathtools}
\usepackage[shortlabels]{enumitem}
\usepackage{xcolor}
\usepackage[colorlinks=true,linkcolor=blue!50!black,citecolor=blue!50!black,urlcolor=blue!50!black]{hyperref}

\newtheorem{theorem}{Theorem}[section]
\newtheorem{lemma}[theorem]{Lemma}

\theoremstyle{definition}

\theoremstyle{remark}
\newtheorem{remark}[theorem]{Remark}

\newcommand{\Z}{\mathbb{Z}}
\newcommand{\R}{\mathbb{R}}
\newcommand{\Q}{\mathbb{Q}}
\newcommand{\sm}{\setminus}

\newcommand{\ceil}[1]{\lceil#1\rceil}
\newcommand{\st}{\,:\,}
\newcommand{\one}{\mathbf{1}}
\newcommand{\Lnat}{\mathrm{L}^{\natural}}
\DeclareMathOperator{\dom}{dom}

\title{Mixed-integer dual network problems are solvable\\ in strongly polynomial time }
\author{Giacomo Zambelli\thanks{London School of Economics and Political Science, Houghton Street, London WC2A 2AE, United Kingdom. E-mail: g.zambelli@lse.ac.uk}}
\date{November 10, 2026}

\begin{document}
\maketitle

\abstract{In this paper we show that mixed-integer linear programming is solvable in strongly polynomial time if the constraint matrix is totally unimodular with at most two nonzero entries in every row. This is done by reducing the problem to an instance of minimization of an $L$-convex function. }

\section{Introduction}\label{sec:intro}

Let $\mathcal{T}$ be the class of totally unimodular matrices in which every row has at most two nonzero entries, and $[n]:=\{1,\ldots,n\}$. The following is the main result of the paper. 
\begin{theorem}\label{thm:main}
There exists and algorithm that, given an $m\times n$ matrix in $\mathcal T$, $b\in Q^m$, $c\in \Q^n$ and $I\subseteq [n]$, solves the mixed-integer linear programming problem
\begin{equation}\label{eq:main}
\min\left\{c^\top x\st  Ax\ge b,\, x_i\in\Z\ (i\in I)\right\}
\end{equation}
in time $O(n \log n (m + n \log n)k^4 \log^2 k)$, where $k=|I|$.
\end{theorem}

The cases $I=\emptyset$ and $I=[n]$ can be cast as the dual of a minimum-cost flow problem, where in the case $I=[n]$ we only need to replace $b_e$ with $\ceil{b_e}$ ($e\in E$). Similarly, the problem reduces to a linear programme if $b\in\Z^m$. To the best of the author's knowledge, a polynomial-time algorithm was not previously know for the general case. The algorithm is based on a reduction to the minimization of an $L$-convex function. $L$-convex functions were introduced by Murota \cite{Murota98} as a discrete notion of convex functions and as an extension of the concept of submodular functions to the integer lattice. We refer the reader to Murota's monograph~\cite{Mur03} for an extensive treatment on the topic.

Problem~\eqref{eq:main} was introduced by Conforti, Di Summa, Eisenbrand and Wolsey
\cite{CDEW09}. Their primary interest was characterizing the convex hull of solutions, 
as a generalization of several mixed-integer
sets considered in the literature \cite{ConDiSWol,GuPo01,MilWol03b,PocWol94,VV}, usually as relaxations of lot-sizing problems \cite{PocWol06}. Specifically, 
Conforti et al.~\cite{CDEW09} gave an extended formulation for the convex hull whose size, while polynomial in $m$ and $n$, may depend exponentially on the encoding size of $b$. In particular, their extended formulation implies that,  whenever the least common denominator of the numbers $b_e$ ($e\in E$) is bounded by a constant, problem~\eqref{eq:main} can be solved in polynomial time.  It should be pointed out that the feasibility problem for the natural ``dual'' counterpart of ~\eqref{eq:main} is NP-hard; indeed, Conforti et al.~\cite{CDEW09} proved that finding a solution to $Cx=d$, $x\ge 0$, $x_i\in \Z$ ($i\in I$) is NP-complete in general if $C$ totally unimodular with two nonzeros per column, even if $d$ is half-integral.

In the next section we define $L$-convex functions and summarize the necessary background on $L$-convex function minimization. In Section~\ref{sec:reduction} we first observe that we can reduce~\eqref{eq:main} to the case where $A$ has exactly a $+1$ and a $-1$ in every row, and then prove the $L$-convexity of the function defined by $g(p)=\{c^\top x\st  Ax\ge b,\, x_i=p_i\ (i\in I)\}$ for  $p\in\Z^I$. The running time bound in Theorem~\ref{thm:main} will follow from an algorithm of Kolmogorov and Shioura~\cite{kolmogorov2009} for $L$-convex function minimization, Lee et al.'s algorithm for submodular function minimization~\cite{lee2015faster}, an algorithm of Orlin for minimum-cost flow~\cite{Orlin93}, and a proximity result of Paat et al.~\cite{paat2020distances}.

\paragraph{Notations.}
Given a finite set $N$, we will denote by $\one_N$ the all-ones vector in $\Z^N$. When there is no ambiguity, we will often simply write $\one$. Given $p,q\in\R^N$, we denote by $p\vee q$ and $p\wedge q$ the vectors in $\R^N$ defined as the componentwise maximum and minimum, respectively, of $p$ and $q$. Given a function $g:\Z^N\to\R\cup\{+\infty\}$, we denote its effective domain by $\dom g:=\{p\in \Z^N\st g(p)<+\infty\}$. For a vector $x\in\R^N$ and $S\subseteq N$, we denote by $x_S$ the restriction of $x$ to the indices in $S$. Similarly, given a matrix $A$ and a subset $S$ of the indices of the columns of $A$, we will denote by $A_S$ the submatrix of $A$ defined by the columns of $A$ indexed by $S$. 

\section{$L$-convex functions}

Given a finite ground set $N$, a function $g:\Z^N\to\R\cup\{+\infty\}$ is {\em $L$-convex} if it satisfies the following conditions
\begin{align}
g(p)+g(q) \ge\ g\big(p\vee q\big)+g\big(p\wedge q\big) &\quad \forall p,q\in \Z^N,\label{eq:L_submodular}\\
\exists r\in\R \mbox{ such that } g(p+\one)=g(p)+r &\quad \forall p\in\Z^N.\label{eq:L-linear}
\end{align}

Property \eqref{eq:L_submodular} is {\em  submodularity} on the integer lattice, property~\eqref{eq:L-linear} is {\em linearity} in the direction of $\one$.  We refer the reader to Chapter 7 of Murota's book~\cite{Mur03} for a treatment of the theory of $L$-convex functions.

The {\em $L$-convex minimization problem} consists of computing $p^*\in\Z^N$ such that $g(p^*)\le g(p)$ for all $p\in\Z^N$. The problem is interesting only when $g$ is constant in the direction of $\one$,
\[g(p+\one)=g(p) \quad \forall p\in\Z^N,
\]
otherwise $g$ does not have a minimum

Murota's {\em scaling steepest-descent} algorithm~\cite{murota2003steepest} (see also~\cite[Section 10.3.2]{Mur03} ) can be used to find a minimizer to a given $L$-convex function, provided that the effective domain is bounded and that a starting point is provided. 

The algorithm requires, at every iteration, the minimization of a submodular function defined by $g$ and by the current point $p\in\dom g$, namely
\[
\rho_p(S)=g(p+\chi_S)-g(p).
\]
The {\em submodular function minimization problem} can be solved in strongly polynomial time, that is, the number of arithmetic operations and oracle calls to $g$ is polynomially bounded by $n=|N|$~\cite{IFF01,Sch00}.  Currently, the best upper bound on the complexity of submodular function minimization is $O(n^3 \log^2 n\cdot EO + n^4 \log^{O(1)} n)$, due to Lee et al.~\cite{lee2015faster} (see also~\cite{dadush2021geometric}), where $EO$ denotes the time required for a single call to the evaluation oracle of the submodular function.

The best running time bound for an algorithm to minimize $L$-convex functions is due to Kolmogorov and Shioura\footnote{Kolmogorov and Shioura's algorithm solves the minimization problem for a larger class of functions, namely $\Lnat$-convex functions, so their result is stated for this larger class in~\cite{kolmogorov2009}.}~\cite{kolmogorov2009}. We denote by $\tau(n)$ the running time required to compute the minimum of a submodular function on $n$ elements.
\begin{theorem}[Kolmogorov and Shioura~{\cite[Theorem 2.11]{kolmogorov2009}}] \label{thm:Lnat-minimization}
Let $g:\Z^N\to\R\cup\{+\infty\}$ be $L$-convex and $K=\max \{\|p-q\|_\infty\st p,q\in \dom g\}$. If $K<+\infty$ and a starting point $p\in\dom g$ is provided, there is an algorithm that finds an optimum $p^*$ to $\min_{p\in\Z^N} g(p)$ in time $O(n\max\{1,\log(K/2n)\}\tau(n))$ where $n=|N|$.
\end{theorem}

\section{Proof of Theorem~\ref{thm:main}}\label{sec:reduction}

In this section we show that~\eqref{eq:main} can be reduced to the problem of minimizing a suitable $L$-convex function. First, we will observe that~\eqref{eq:main} can be written in the following form. Let $G=(V,E)$ be a directed graph, $b\in \Q^E$, $c\in\Q^V$, $I\subseteq V$ such that $c^\top \one=0$. We consider the following problem.

\begin{equation}\label{eq:reduced}
\min\left\{c^\top x\st  x_v-x_u\ge b_{uv}\ (u,v)\in E,\,  x_I\in\Z^I\right\}.
\end{equation}
We will show that every problem~\eqref{eq:main} can be reduced to a problem of the form~\eqref{eq:reduced}, and then provide a polynomial-time algorithm to solve \eqref{eq:reduced} by a reduction to \emph{$L$-convex} function minimization. 

\begin{remark}\label{rmk:lineality}
The lineality space of the feasible region of the LP relaxation of~\eqref{eq:reduced} contains the line generated by $\one$, since this is a solution to $x_v-x_u=0$ ($(u,v)\in E$). Therefore, in order for~\eqref{eq:reduced} to have bounded optimal value, the hypothesis that $c^\top \one=0$ is necessary.
\end{remark}
\medskip

We now explain the reduction from~\eqref{eq:main} to~\eqref{eq:reduced}. Consider $A$, $b$, $c$, $I$ defined as in the statement of Theorem~\ref{thm:main}. Since $A$ is totally unimodular with at most two non-zero entries per row, by a theorem of Heller and Tompkins~\cite{heller1956}, $[n]$ can be partitioned into sets $B$ and $R$ such that  the sum of the columns of $A$ indexed by $B$ minus the sum of the columns indexed by $R$ is a $0,\pm 1$ vector, that is, $\bar a:=A_B\one_B-A_R\one_R\in\{0,+1,-1\}^m$. Let $\gamma:=c^\top_B\one_B-c^\top_R\one_R$, and define $\bar A=( A_B,-A_R,-\bar a)$ and $\bar c^\top=(c_B^\top, -c_R^\top,-\gamma)$. Note that~\eqref{eq:main} is equivalent to  $\min\{\bar c^\top z\st \bar Az\ge b,\, z_I\in\Z^I,\, z_{n+1} = 0\}$ via the change of variables $x_i=z_i$ ($i \in B$), $x_i=-z_i$ ($i \in R$). Let $V=[n+1]$ and $J=I\cup\{n+1\}$, and note that $\bar A\one_V=0$ and $\bar c^\top\one_V=0$. This implies that $\min\{\bar c^\top z\st \bar Az\ge b,\, z_I\in\Z^I,\, z_{n+1} = 0\}$ is equivalent to $\min\{\bar c^\top z\st \bar Az\ge b,\,z_{J}\in\Z^{J}\}$. Indeed, given any optimal solution $z^*$ to the latter problem, all points of the form $z^*+\alpha\one_V$ ($\alpha\in\Z)$ are optimal as well, hence by choosing $\alpha=-z^*_{n+1}$ we can ensure that $z^*_{n+1}=0$. Finally, since $\bar A$ has exactly a $+1$ and a $-1$ entry in every row, it is the transpose of the incidence matrix of a directed graph $G=(V,E)$, and $\min\{\bar c^\top z\st \bar Az\ge b,\, z_J\in\Z^J\}$ is a problem of the form~\eqref{eq:reduced}.
\bigskip

In the reminder of this section, we consider problem~\eqref{eq:reduced}, and denote by $X\subseteq\R^V$ its feasible region. 

\begin{lemma}
\label{lemma:preserve-X} Let $x,y\in \R^V$. Let $x'=x\vee y$ and $y'=x\wedge y$.
The following hold.
\begin{enumerate}[(i)]
  \item $x+y=x'+y'$,
  \item If $x,y\in X$, then $x',y'\in X$. 
\end{enumerate}  
\end{lemma}
\begin{proof}
Part (i) follows immediately from the fact that $\{x_v,y_v\}=\{x_v\vee y_v,x_v\wedge y_v\}$ for all $v\in V$. For (ii), assume that $x,y\in X$. Clearly $x'_I,y'_I\in\Z^I$.  Let $(u,v)\in E$. We have $x'_v-x'_u\ge \max\{x_v-x_u,y_v-y_u\}\ge b_{uv}$ and  $y'_v-y'_u\ge \max\{x_v-x_u,y_v-y_u\}\ge b_{uv}$
\end{proof}

If~\eqref{eq:reduced} is feasible, then it has an optimum if and only if its LP relaxation has an optimum. The LP relaxation of~\eqref{eq:reduced} is the dual of a minimum-cost flow problem  an can be solved by an algorithm of Orlin~\cite{Orlin93} in time $O(n \log n (m + n \log n))$. We can therefore ensure that~\eqref{eq:reduced} is not unbounded. Under such hypothesis, we define a function $g\st \Z^I\to \R\cup\{+\infty\}$ by
\begin{equation}\label{eq:g-function}
g(p)=\min\left\{c^\top x\st x_v-x_u\ge b_{uv}\ (u,v)\in E,\, x_I=p\right\}\qquad  \forall p\in\Z^I,
\end{equation}
where $g(p)=+\infty$ if the minimization problem in \eqref{eq:g-function} is infeasible. Note that, by the assumption that the LP relaxation of~\eqref{eq:reduced} is not unbounded, $g(p)$ has finite value for every $p$ for which \eqref{eq:g-function} is feasible.

\begin{lemma}\label{lemma:g-L-natural}
Assume that problem~\eqref{eq:reduced} admits an optimum, and denote by $\omega^*$ the optimum value. The function $g\st \Z^I\to \R_+$ defined in~\eqref{eq:g-function} is $L$-convex, constant in the direction of $\one_I$, and $\dom g\neq\emptyset$. Furthermore, $\min_{p\in\Z^I} g(p)=\omega^*$. 
\end{lemma}
\begin{proof}
Since problem~\eqref{eq:reduced} admits an optimum, then there is at least a point for which~\eqref{eq:g-function} is feasible, hence  $\dom g\neq \emptyset$. By definition, $f(p)\ge \omega^*$ for every $p\in\Z^I$ and, given an optimum $x^*$ for~\eqref{eq:reduced}, the point $p^*=x^*_I$ satisfies $f(p^*)=\omega^*$. It follows that  $\min_{p\in\Z^I} f(p)=\omega^*$. 

It remains to show that $g$ is $L$-convex. Observe that, since $c^\top\one_V=0$, by Remark~\ref{rmk:lineality} $g(p+\one_I)=g(p)$ for all $p\in\Z^I$, hence $g$ is constant in the direction of $\one_I$, so in particular it satisfies~\eqref{eq:L-linear} for $r=0$.
We need to verify that $g$ is submodular over $\Z^I$, that is, it satisfies~\eqref{eq:L_submodular}. Let  $p,q\in z^I$. We may assume $g(p), g(q)<+\infty$, otherwise ~\eqref{eq:L_submodular} is trivially satisfied. It follows that there exist $x,y\in X$ such that $x_I=p$, $y_I=q$, $g(p)=c^\top x$, $g(q)=c^\top y$. Let $x'=x\vee y$ and $y'=x\wedge y$. Note that $p\vee q=x'_I$ and $p\wedge q=y'_I$. By Lemma~\ref{lemma:preserve-X}, $x',y'\in X$.  It follows that 
\[g(p\vee q)+g(p\wedge q)\le c^\top x'+ c^\top y'= c^\top x+ c^\top y=g(p)+g(q),\]
where the first equality follows because $x+y=x'+y'$ by Lemma~\ref{lemma:preserve-X}(i). 
\end{proof}

\begin{lemma}\label{lemma:Orlin}
For every $p\in\Z^I$, the value of $g(p)$ can be computed in time $O(n \log n (m + n \log n))$, where  $n=|V|$ and $m=|E|$.
\end{lemma}
\begin{proof}
The problem~$\min\{c^\top x\st x\in X, x_I=p\}$ can be cast as the dual of a minimum-cost flow problem with $O(n)$ nodes and $O(m)$ arcs. Specifically, for every $e=(u,v)\in E$, let $\bar b_e=b_e$ if $u,v\in V\sm I$, let $\bar b_e=b_e+p_u$ if $u\in I ,v\in V\sm I$, and $\bar b_e=b_e-p_v$ if $u\in V\sm I ,v\in I$. Let $G'=(V',E')$ be the directed graph obtained by removing all arcs with both endnodes in $I$ and contracting $I$ into a single node $v^I$. Finally, define $\bar c_v=c_v$ ($v\in V\sm I$) and $\bar c_{v^I}=\sum_{v\in I}c_v$. Observe that the problem $g(p)=\min\left\{c^\top x\st x_v-x_u\ge b_{uv}\ (u,v)\in E,\, x_I=p\right\}$ is equivalent to $\min\left\{\bar c^\top x\st x_v-x_u\ge \bar b_{uv}\ (u,v)\in E',\, x_{v^I}=0\right\}$. Noting that, by construction, $\bar c^\top \one_{V'}=c^\top \one_V=0$, the latter problem is equivalent to  $\min\left\{\bar c^\top x\st x_v-x_u\ge \bar b_{uv}\ (u,v)\in E'\right\}$, since we can ensure that an optimal solution with $x_{v^I}=0$ exists whenever an optimum exists. The latter problem is the dual of a minimum-cost flow problem, hence it can be solved by Orlin's algorithm~\cite{Orlin93} in time $O(n \log n (m + n \log n))$.
\end{proof}

By the above lemma, solving~\eqref{eq:reduced} is equivalent to solving  $\min_{p\in\Z^I} g(p)$. To apply the $L$-convex function minimization algorithm of Theorem~\ref{thm:Lnat-minimization}, we need the domain of $g$ to have bounded $\ell_\infty$ diameter. We will achieve this by restricting the domain of $g$ to a box that is guaranteed to have an optimal solution. For this, we will appeal to the following proximity result.

\begin{theorem}[Paat et. al \cite{paat2020distances}] \label{thm:Paat} Let $C$ be an integer $m\times n$ matrix, $I\subseteq [n]$, $d\in\Q^m$, $w\in\R^n$. For every optimal solution $\bar x$ to $\max\{w^\top x\st Cx\le d\}$ there exists an optimal solution $x^*$ to $\max\{w^\top x\st Cx\le d,\, x_I\in\Z^I\}$ such that $\|\bar x-x^*\|_\infty\le |I|\Delta$, where $\Delta$ is the maximum of the absolute values of the determinants of the square submatrices of $C$.
\end{theorem}
The above improves the classic  $n\Delta$ bound of Cook et. al \cite{CGST86} by making it dependent only on the number of integer variables\footnote{It should be noted that the result as stated in~\cite{paat2020distances} declares $d\in\Z^m$, but the authors point out that indeed their proof implies the results for any $d\in \Q^m$}.

\begin{theorem}
\label{thm:phase_2} Assume that the LP relaxation of~\eqref{eq:reduced} admits an optimum  $\bar x$, and let $\mathcal{D}:=\dom(g)\cap \{p\in\Z^I\st \bar x_i-|I|\le p_i\le \bar x_i+|I|\ (i\in I)\}$. Let $n=|V|$, $m=|E|$, $k=|I|$.
\begin{enumerate}[(i)]
  \item There is an algorithm that, given $p\in \mathcal{D}$, computes an optimum $x^*$ of~\eqref{eq:reduced} with $x^*_I\in\mathcal{D}$ in time $O(n \log n (m + n \log n)k^4 \log^2 k)$.
  \item There is an algorithm that, in time $O(n \log n (m + n \log n)k^4 \log^2 k)$,  computes a point $p\in \mathcal{D}$ or correctly determines that~\eqref{eq:reduced} is infeasible.
\end{enumerate}

\end{theorem}
\begin{proof}
(i) Let us define the function $f\st \Z^I\to \R$ to be the restriction of $g$ to $\mathcal{D}$, that is, $f(p)=g(p)$ for $p\in \mathcal{D}$, $f(p)=+\infty$ for $p\in \Z^I\sm \mathcal{D}$. By Theorem~\ref{thm:Paat}, since the constraint matrix of~\eqref{eq:reduced} is totally unimodular, the problem admits an optimum $x^*$ such that $\|\bar x-x^*\|_\infty\leq k$, so in particular $x^*_I\in\mathcal{D}$, since $x^*_I$ is integral. This implies that $\dom f$ contains an optimal solution to the problem $\min_{p\in \Z^I} g(p)$. In particular, $\dom f\neq \emptyset$, hence Murota~\cite[Theorem 7.11(2)]{Mur03} implies that $f$ is $L$-convex. By Theorem~\ref{thm:Lnat-minimization}, a minimizer $p^*$ for $f$ can be computed in time $O(k\tau(k))$, because $\max \{\|p-q\|_\infty\st p,q\in \dom f\}\le 2k$. By our choice of $\mathcal{D}$, $p^*$ is also a minimizer of $g$, hence a solution $x^*$ to~\eqref{eq:reduced} can be computed by solving the LP~\eqref{eq:g-function} for $p=p^*$.    

For the overall running time, recall that a single submodular function minimization problem over the ground set $I$ can be solved in time $O(k^3 \log^2 k\cdot EO + k^4 \log^{O(1)} k)$ by the algorithm of Lee et al.~\cite{lee2015faster}, and every evaluation oracle call can  be computed in time $O(n \log n (m + n \log n))$ as in Lemma~\ref{lemma:Orlin}. Since $k\le n$, $O(k^4 \log^{O(1)} k)$ is dominated by $O(n \log n (m + n \log n)k^4 \log^2 k)$, hence providing the stated  time-complexity bound.

\noindent (ii) Define a new digraph $\bar G=(\bar V,\bar E)$ by adding to $G$ a node $v'$ for every $v\in I$ and the arc $(v,v')$, so $\bar V=V\cup J$ where $J=\{v'\st v\in I\}$ and $\bar E=E\cup F$ where $F=\{(v,v')\st v\in I\}$. Let $\bar b\in\Q^{\bar E}$ be defined by $\bar b_E=b$ and $\bar b_F=0$, and let $\bar c\in\Q^{\bar V}$ be defined by $\bar c_v=1$ for $v\in J$, $\bar c_{v}=-1$ for  $v\in I$, and $c_v=0$ otherwise. Note that $\bar c^\top y=\sum_{v\in I} y_{v'}-y_{v}$. Consider the instance of~\eqref{eq:reduced} defined by $\bar G$, $\bar b$, $\bar c$ and $J$, that is 
\begin{equation}\label{eq:phase2}
\min\{\bar c^\top y\st  y_v-y_u\ge \bar b_{uv}\ (u,v)\in \bar E,\, y_J\in\Z^J\}.
\end{equation}

Note that, by construction, $\bar c^\top y\ge 0$ for every feasible solution of the LP relaxation of the latter problem. If we define $\bar y\in\R^{\bar V}$ by $\bar y_V=\bar x$, $\bar y_{v'}=\bar x_v$ for all $v\in I$, $\bar y$ is feasible for the LP relaxation of \eqref{eq:phase2} and it has value $0$, hence it is an optimal solution for the LP relaxation. Furthermore, note that the point $\tilde{y}\in\R^{\bar v}$ defined by $\tilde y_V=\bar y_V$, $\tilde y_J=\ceil{\bar y_J}$ is feasible for~\eqref{eq:phase2}. 

The algorithm in (i) computes an optimal solution $y^*$ to~\eqref{eq:phase2} such that $\bar y_{v'}-k\le y^*_{v'}\le\bar y_{v'}+k$ for all $v\in I$ in time $O(n \log n (m + n \log n)k^4 \log^2 k)$. Note that $\bar c^\top y^*\ge 0$, and $\bar c^\top y^*=0$ if and only if it satisfies
\[y_v-y_u\ge b_{uv}\ (u,v)\in E,\, y_I=y_J,\, y_J\in\Z^J.\]
If so, it follows that $p=y^*_I$ is in $\mathcal{D}$. If $\bar c^\top y^*>0$, it follows that~\eqref{eq:reduced} is infeasible. 
\end{proof}

\section{Discussion}

Ultimately, the algorithm presented reduces problem~\eqref{eq:main} to a sequence of minimum-cost flow computations. A natural question is whether a more direct reduction to minimum-cost flow exists. In particular, the extended formulation of Conforti et al.~\cite{CDEW09} reduces the problem to a single minimum-cost flow computation on a graph with $O(nD)$ nodes and $O(mD)$ arcs, where $D$ is the least common denominator of $b_e$ ($e\in E$). It is an open question whether a reduction to a minimum-cost flow on a network of polynomial size exists.  

Another outstanding open problem is to provide an inequality description in the original space of variables of the convex hull of mixed-integer solutions. This is known only for the case that $b$ is half-integral~\cite{CGZ07} or for specific topologies of the underlying graph of the constraint matrix $A$ \cite{ConDiSWol,CWZ10,GuPo01,MilWol03b,PocWol94,VV}. Interestingly, a connection with submodularity in the context of characterizing the mixed-integer hull  had previously been observed by K{\i}l{\i}n\c{c}-Karzan et al.~\cite{KilKucLee22} in the special case of the so-called {\em mixing set} (see~\cite{GuPo01}). Given the connection to $L$-convexity presented here for the optimization problem, it is natural to ask if this has implications for the description of the mixed-integer hull.

\subsection{Statement of AI use}

The paper originates from the following insight obtained using Anthropic Opus 5.5, reducing a special case of~\eqref{eq:main} to submodular function minimization. Specifically, given a bipartite graph $G$, and assuming $0\le b_e\le 1$, consider the following mixed-integer bipartite vertex cover problem $\min\{c^\top x\st x_u+x_w\ge b_{uw}\ (uw\in E),\, x_U\ge 0,\, x_I\in\Z^I\}$. The problem has an optimum if and only if $c\ge 0$, hence there exists an optimal solution with $x_I\in\{0,1\}^I$. If $U,W$ is the bipartition of $G$, and we perform a change of variables, replacing $x_w$ with $-x_w$ for all $w\in W\sm I$, $x_w$ with $1-x_w$ for $w\in W\cap I$, we obtain an equivalent problem of the form $\min\{\bar c^\top x\st x_u-x_w\ge \bar b_{uw}\ (uw\in E),\, x_U\ge 0,\, x_{W\sm I}\le 0,\, x_I\in\{0,1\}^I\}$ (where $u$ and $w$ denote nodes in $U$ and $W$, respectively). For $S\subseteq I$, let $f(S)$ be the optimal value of the latter problem once we fix $x_v=1$ if $v\in S$, $x_v=0$ if $v\in I\sm S$. Anthropic Opus 5.5 observed that the function $f\st 2^I\to\R$ is submodular, therefore reducing this special case to a single instance of submodular function minimization. 

The generalization of the reduction through the use of $L$-convex functions, as well as strong polynomiality through the use of proximity, was the work of the author alone. The paper was entirely written by the author.

\bibliographystyle{plain}
\bibliography{references}

@string{mor      = "Math. Oper. Res."}

@string{mp       = "Math. Program."}

@string{jacm     = "J. ACM"}

@string{discopt  = "Discrete Optim."}

@string{siamopt  = "SIAM J. Optim."}

@string{orjnl    = "Oper. Res."}

@string{jctb     = "J. Combin. Theory Ser. B"}

@article{CDEW09,
  author  = {Conforti, M. and {Di Summa}, M. and Eisenbrand, F. and Wolsey, L. A.},
  title   = {Network formulations of mixed-integer programs},
  journal = mor,
  volume  = {34},
  pages   = {194--209},
  year    = {2009}
}

@inproceedings{ConDiSWol,
  author    = {Conforti, M. and {Di Summa}, M. and Wolsey, L. A.},
  title     = {The intersection of continuous mixing polyhedra and the continuous mixing polyhedron with flows},
  booktitle = {Integer Programming and Combinatorial Optimization (IPCO 2007)},
  editor    = {Fischetti, M. and Williamson, D. P.},
  series    = {Lecture Notes in Computer Science},
  volume    = {4513},
  pages     = {352--366},
  publisher = {Springer},
  address   = {Berlin, Heidelberg},
  year      = {2007}
}

@inproceedings{CGZ07,
  author    = {Conforti, M. and Gerards, A. M. H. and Zambelli, G.},
  title     = {Mixed-integer vertex covers on bipartite graphs},
  booktitle = {Integer Programming and Combinatorial Optimization (IPCO 2007)},
  editor    = {Fischetti, M. and Williamson, D. P.},
  series    = {Lecture Notes in Computer Science},
  volume    = {4513},
  pages     = {324--336},
  publisher = {Springer},
  address   = {Berlin, Heidelberg},
  year      = {2007}
}

@article{CWZ10,
  author  = {Conforti, M. and Wolsey, L. A. and Zambelli, G.},
  title   = {Projecting an extended formulation for mixed-integer covers on bipartite graphs},
  journal = mor,
  volume  = {35},
  pages   = {603--623},
  year    = {2010}
}

@article{CGST86,
  author  = {Cook, W. and Gerards, A. M. H. and Schrijver, A. and Tardos, {\'E}.},
  title   = {Sensitivity theorems in integer linear programming},
  journal = mp,
  volume  = {34},
  pages   = {251--264},
  year    = {1986}
}

@article{dadush2021geometric,
  author  = {Dadush, D. and V{\'e}gh, L. A. and Zambelli, G.},
  title   = {Geometric rescaling algorithms for submodular function minimization},
  journal = mor,
  volume  = {46},
  pages   = {1081--1108},
  year    = {2021}
}

@article{GuPo01,
  author  = {G{\"u}nl{\"u}k, O. and Pochet, Y.},
  title   = {Mixing mixed integer inequalities},
  journal = mp,
  volume  = {90},
  pages   = {429--457},
  year    = {2001}
}

@incollection{heller1956,
  author    = {Heller, I. and Tompkins, C. B.},
  title     = {An extension of a theorem of {D}antzig's},
  booktitle = {Linear Inequalities and Related Systems},
  editor    = {Kuhn, H. W. and Tucker, A. W.},
  series    = {Annals of Mathematics Studies},
  number    = {38},
  pages     = {247--254},
  publisher = {Princeton University Press},
  address   = {Princeton},
  year      = {1956}
}

@article{IFF01,
  author  = {Iwata, S. and Fleischer, L. and Fujishige, S.},
  title   = {A combinatorial strongly polynomial algorithm for minimizing submodular functions},
  journal = jacm,
  volume  = {48},
  pages   = {761--777},
  year    = {2001}
}

@article{KilKucLee22,
  author  = {K{\i}l{\i}n{\c{c}}-Karzan, F. and K{\"u}{\c{c}}{\"u}kyavuz, S. and Lee, D.},
  title   = {Joint chance-constrained programs and the intersection of mixing sets through a submodularity lens},
  journal = mp,
  volume  = {195},
  pages   = {283--326},
  year    = {2022}
}

@article{kolmogorov2009,
  author  = {Kolmogorov, V. and Shioura, A.},
  title   = {New algorithms for convex cost tension problem with application to computer vision},
  journal = discopt,
  volume  = {6},
  pages   = {378--393},
  year    = {2009}
}

@inproceedings{lee2015faster,
  author       = {Lee, Y. T. and Sidford, A. and Wong, S. C.},
  title        = {A faster cutting plane method and its implications for combinatorial and convex optimization},
  booktitle    = {Proceedings of the 56th Annual {IEEE} Symposium on Foundations of Computer Science ({FOCS})},
  pages        = {1049--1065},
  publisher    = {IEEE},
  year         = {2015}
}

@article{MilWol03b,
  author  = {Miller, A. and Wolsey, L. A.},
  title   = {Tight formulations for some simple {MIP}s and convex objective integer programs},
  journal = {Math. Program. Ser. B},
  volume  = {98},
  pages   = {73--88},
  year    = {2003}
}

@article{Murota98,
  author  = {Murota, K.},
  title   = {Discrete convex analysis},
  journal = mp,
  volume  = {83},
  pages   = {313--371},
  year    = {1998}
}

@book{Mur03,
  author    = {Murota, K.},
  title     = {{Discrete Convex Analysis}},
  series    = {SIAM Monographs on Discrete Mathematics and Applications},
  publisher = {SIAM},
  address   = {Philadelphia},
  year      = {2003}
}

@article{murota2003steepest,
  author  = {Murota, K.},
  title   = {On steepest descent algorithms for discrete convex functions},
  journal = siamopt,
  volume  = {14},
  pages   = {699--707},
  year    = {2004}
}

@article{Orlin93,
  author  = {Orlin, J. B.},
  title   = {A faster strongly polynomial minimum cost flow algorithm},
  journal = orjnl,
  volume  = {41},
  pages   = {338--350},
  year    = {1993}
}

@article{paat2020distances,
  author  = {Paat, J. and Weismantel, R. and Weltge, S.},
  title   = {Distances between optimal solutions of mixed-integer programs},
  journal = mp,
  volume  = {179},
  pages   = {455--468},
  year    = {2020}
}

@article{PocWol94,
  author  = {Pochet, Y. and Wolsey, L. A.},
  title   = {Polyhedra for lot-sizing with {Wagner--Whitin} costs},
  journal = mp,
  volume  = {67},
  pages   = {297--323},
  year    = {1994}
}

@book{PocWol06,
  author    = {Pochet, Y. and Wolsey, L. A.},
  title     = {{Production Planning by Mixed Integer Programming}},
  series    = {Springer Series in Operations Research and Financial Engineering},
  publisher = {Springer},
  address   = {New York},
  year      = {2006}
}

@article{Sch00,
  author  = {Schrijver, A.},
  title   = {A combinatorial algorithm minimizing submodular functions in strongly polynomial time},
  journal = jctb,
  volume  = {80},
  pages   = {346--355},
  year    = {2000}
}

@article{VV,
  author  = {{Van Vyve}, M.},
  title   = {The continuous mixing polyhedron},
  journal = mor,
  volume  = {30},
  pages   = {441--452},
  year    = {2005}
}

\end{document}